\documentclass[english,a4paper,12pt]{smfart}

\usepackage[french]{babel}
\usepackage{amssymb}
\usepackage{latexsym}
\usepackage{bm}

\numberwithin{equation}{section}

\def\di{\displaystyle}

\def\un{\underline}

\def\bn{\mbox{\bf n}}
\def\ba{\mbox{\bf a}}
\def\bb{\mbox{\bf b}}
\def\bbm{\mbox{\bf m}}

\def\Pren{\mbox{\rm Pren}}
\def\Prenb{\mbox{\rm Pren}^{\bullet}}
\def\Mb{\mbox{\rm M}^{\bullet}}
\def\M{\mbox{\rm M}}
\def\Nb{\mbox{\rm N}^{\bullet}}
\def\N{\mbox{\rm N}}

\def\Flin{\mbox{\rm F}_{\rm lin}}
\def\FSem{\mbox{\rm F}_{\rm Sem}}
\def\Fsem{\mbox{\rm F}_{\rm sem}}
\def\FTrem{\mbox{\rm F}_{\rm Trem}}
\def\FPoin{\mbox{\rm F}_{\rm Poin}}

\def\FDulac{\mbox{\rm F}_{\rm Dulac}}

\def\1b{\mbox{\rm 1}^{\bullet}}
\def\un{\mbox{\rm 1}}
\def\I{\mbox{\rm I}}
\def\Ib{\mbox{\rm I}^{\bullet}}
\def\Bb{\mbox{\rm B}_{\bullet}}
\def\Db{\mbox{\rm D}_{\bullet}}

\def\Conjb{\mbox{\rm C}^{\bullet}}
\def\Bbb{\mbox{\rm\bf B}}
\def\Thetab{\Theta^{\bullet}}
\def\Thetaib{{\bm{-\!\!\!-}}\!\!\!\!\!\!\Theta}
\def\Mib{{\bm{-\!\!\!-}}\!\!\!\!\!\!\mbox{\rm M}^{\bullet}}
\def\dem{\mbox{\rm dem}}
\def\demb{\mbox{\rm dem}^{\bullet}}
\def\Dem{\mbox{\rm Dem}}
\def\Demb{\mbox{\rm Dem}^{\bullet}}

\def\denb{\mbox{\rm den}^{\bullet}}

\def\Denb{\mbox{\rm Den}^{\bullet}}

\def\sumb{\sum_{\bullet}}
\def\edelta{\mbox{\rm e}^{\Delta}}

\def\sem{\mbox{\rm sem}}
\def\semb{\mbox{\rm sem}^{\bullet}}
\def\Sem{\mbox{\rm Sem}}
\def\Semb{\mbox{\rm Sem}^{\bullet}}

\def\poin{\mbox{\rm poin}}
\def\poinb{\mbox{\rm poin}^{\bullet}}
\def\Poin{\mbox{\rm Poin}}
\def\Poinb{\mbox{\rm Poin}^{\bullet}}

\def\Log{\mbox{\rm Log}}
\def\Exp{\mbox{\rm Exp}}
\def\AutSimp{\mbox{\rm Simp}}

\def\tremb{\mbox{\rm trem}^{\bullet}}
\def\trem{\mbox{\rm trem}}
\def\Trem{\mbox{\rm Trem}}
\def\Tremb{\mbox{\rm Trem}^{\bullet}}

\def\dulacb{\mbox{\rm dulac}^{\bullet}}
\def\dulac{\mbox{\rm dulac}}
\def\Dulac{\mbox{\rm Dulac}}
\def\Dulacb{\mbox{\rm Dulac}^{\bullet}}

\def\Zib{\mbox{\rm Z}^{\bullet}}
\def\Zi{\mbox{\rm Z}}
\def\limstat{\mbox{\rm limstat}}
\def\indic{\mbox{\bf 1}}

\def\K{\mathbb{K}}

\def\N{\mathbb{N}}
\def\Z{\mathbb{Z}}

\def\C{\mathbb{C}}

\def\cal{\mathcal}
\def\Exp{\mbox{\rm Exp}}

\usepackage{amsmath}
\usepackage{amssymb}
\usepackage{amsthm}
\usepackage{a4wide}

\newtheorem{thm}{Theorem}
\newtheorem{defi}{Definition}
\newtheorem{lem}{Lemma}

\newtheorem{rema}{Remark}

\title[Trimmed and Poincar\'e-Dulac normal form]{About the Trimmed and the Poincar\'e-Dulac normal form of diffeomorphisms}
\alttitle{Trimmed and Poincar\'e-Dulac normal form}
\author{Jacky CRESSON}

\author{Jasmin RAISSY}

\begin{document}
\maketitle
\setcounter{tocdepth}{3}
\baselineskip 6mm

\begin{abstract}
We study two particular continuous prenormal forms as defined by Jean Ecalle and Bruno Vallet for local analytic diffeomorphism of $\C^{\nu}$: the Trimmed form and the
Poincar\'e-Dulac normal form. We first give a self-contain introduction to the mould formalism of Jean Ecalle. We provide a dictionary between moulds and
the classical Lie algebraic formalism using non-commutative formal power series. We then give full proofs and details for results
announced by J. Ecalle and B. Vallet about the Trimmed form of diffeomorphisms. We then discuss a mould approach to the classical Poincar\'e-Dulac normal
form of diffeomorphisms. We discuss the universal character of moulds taking place in normalization problems.
\end{abstract}

\begin{tiny}
\tableofcontents
\end{tiny}

\section{Introduction}

In this paper we study the set of {\it local analytic resonant diffeomorphisms} of $\C^{\nu}$ using the theory of {\it continuous
prenormalization} developped by J. Ecalle (\cite{es},\cite{ev1}). We assume that diffeomorphisms have a diagonalizable
linear part and we work in a chart where the linear part is diagonal. The diffeomorphism is called in this case in
{\it prepared form}. Let $f$ be a diffeomorphism in prepared form. Roughly speaking a {\it prenormal form} of $f$ is
a diffeomorphism $f_{\rm pren}$, conjugated to $f$, of the form $f_{\rm pren}=f_{\rm lin} +f_{\rm rem}$ where $f_{\rm rem}$
is made of resonant terms. A {\it normal form} is a prenormal form containing the minimal number of resonant terms, with
formal invariants as coefficients. Although a normal form can be considered as the simplest prenormal form, it is not
in general possible to compute it. Even if an {\it algorithmic} procedure can be obtained \cite{bai}, its exact shape
is related to the vanishing of certain quantities depending polynomially on the Taylor coefficients of the diffeomorphisms.
This can not be decided by a computer.\\

We look for {\it calculable} prenormal forms, {\it i.e.} prenormal forms which can be obtained using a
procedure which is {\it algorithmic} and {\it implementable}. As an example of such prenormal forms, we study {\it
continuous prenormal forms} as defined by J. Ecalle \cite{es}.\\

We mainly focus on two particular continuous prenormal
forms, one introduced by J. Ecalle and B. Vallet \cite{ev1} called the {\it Trimmed form} and the classical
{\it Poincar\'e-Dulac normal form}. The framework of continuous prenormalization is the {\it mould formalism} developped
by J. Ecalle since 1970. We provide a self-contained introduction to this formalism, omitting some aspects
which will not be used in this paper. We refer to (\cite{cr1},\cite{cr2}) for more details.\\

The Trimmed form is first studied. We give complete proofs for results which are announced by J. Ecalle and B. Vallet
\cite{ev1} with (or without) a sketch of proof. In particular, we give all the details for the computations of the
different moulds associated to the Trimmed form. We also give {\it closed} formulae for these moulds using a different
initial alphabet.\\

The {\it Poincar\'e-Dulac normal form} is then discussed in the mould framework and compared to the Trimmed form. We obtain two universal
moulds $\Poinb$ and $\Dulacb$. These two universal moulds are associated to the Poincar\'e normalization procedure and the Poincar\'e-Dulac
normal form. It seems impossible to obtain such objects using the existing methods of perturbation theory. The mould formalism provides a
direct and algorithmic way to capture the universal features of a normalization procedure.

\section{Diffeomorphisms, automorphisms and continuous prenormalization}

We consider local analytic diffeomorphisms of $\C^{\nu}$ with $0$ as a fixed point and a {\it diagonalizable}
linear part. We work in a given analytic chart where the linear part is assumed to be in diagonal form. In such a case, the
diffeomorphism is called in {\it prepared form} by J. Ecalle.\\

Let $f:\C^{\nu} \rightarrow \C^{\nu}$, $\nu\in \N$ defined by
\begin{equation}
f(x_1 ,\dots ,x_{\nu} )=(e^{\lambda_1} x_1 ,\dots ,e^{\lambda_{\nu}} x_{\nu} ) +h(x_1 ,\dots ,x_{\nu} ) ,
\end{equation}
with $f(0)=0$, and $h=(h_1 ,\dots ,h_{\nu} )$, $h_i \in \C \{ x\}$ for all $i=1,\dots ,\nu$. We denote by $f_{\rm lin}$
the linear part of $f$, {\it i.e.}
$f_{\rm lin} (x_1 ,\dots ,x_{\nu} )=(e^{\lambda_1} x_1 ,\dots ,e^{\lambda_{\nu}} x_{\nu} )$.\\

J. Ecalle looks for the {\it substitution operator} associated to $f$, denoted by $F$ and defined by
\begin{equation}
F:
\left .
\begin{array}{lll}
\C \{ x\} & \rightarrow & \C \{ x\} ,\\
\phi & \mapsto & \phi \circ f ,
\end{array}
\right .
\end{equation}
where $\circ$ is the usual composition of functions.\\

As $f$ is a diffeomorphism, the substitution operator $F$ is an automorphism of $(\C \{ x\} , \cdot )$ where $\cdot$ is the usual
product of functions on $\C \{ x\}$, {\it i.e.} for all $\phi, \psi \in \C \{ x\}$, we have
\begin{equation}
F(\phi \cdot \psi )=F\phi \cdot F\psi ,
\end{equation}
and $F^{-1} (\phi)=\phi\circ f^{-1}$.\\

J. Ecalle uses the following result, which is a direct consequence of the Taylor expansion theorem:

\begin{lem}
Let $f$ be an analytic diffeomorphism of $\C^{\nu}$ in prepared form and $F$ its associated substitution operator. There
exist a decomposition of $F$ as
\begin{equation}
\label{homocomponents}
F=\Flin \left ( \mbox{\rm Id} +\di\sum_{n\in A(F)} B_n \right ) ,
\end{equation}
where $A(F)$ is an infinite set of indices $n\in \Z^{\nu}$, $\Flin$ the substitution operator associated to
$f_{\rm lin}$, and for all $n\in A(F)$, $B_n$ is a homogeneous differential operator of degree $n$, {\it i.e.}
for all $m\in \N^{\nu}$,
\begin{equation}
B_n (x^m )=\beta_{n,m} x^{n+m} ,\ \ \beta_{n,m} \in \C .
\end{equation}
\end{lem}

In the following, we work essentially with the substitution operator $F$. In order to simplify our statements, we call
{\it diffeo(s)} the automorphism $F$ associated to a given diffeomorphism $f$.

\begin{defi}
Let $F$ and $F_{\rm conj}$ be two local analytic diffeos of $\C^{\nu}$. The diffeo $F_{\rm conj}$ is called
conjugated to $F$ if there exists a change of variables $h$ of $C^{\nu}$ such that the associated substitution
operator denoted by $\Theta$ satisfies
\begin{equation}
\label{conjugacy}
F_{\rm conj} =\Theta \cdot F \cdot \Theta^{-1} .
\end{equation}
\end{defi}

The substitution operator $\Theta$ is called the {\it normalizator} in the following. When the change of
variables $h$ is of class formal, $C^k$ or $C^{\omega}$, we speak of a formal, $C^k$ or analytic normalization.\\

\begin{defi}
Let $F$ be an analytic diffeo of $\C^{\nu}$ in prepared form. A prenormal form for $F$, denoted by $F_{\rm pran}$,
is an automorphism of $\C \{ x\}$ conjugated to $F$ such that
\begin{equation}
F_{\rm pran} \cdot \Flin =\Flin \cdot F_{\rm pran} .
\end{equation}
\end{defi}

J. Ecalle has introduced in \cite{es} and extensively studied in \cite{ev1} a very particular class of prenormal forms
called {\it continuous prenormal forms}.

\begin{defi}
Let $F$ be a diffeo of $\C^{\nu}$ in prepared form given by
$$F=\Flin \left ( \mbox{\em Id} +\di\sum_{n\in A(F)} B_n \right ) .$$
A continuous prenormal form $F_{\rm pren}$ is an automorphism of $\C \{ x\}$ of the form
\begin{equation}
\label{normalform}
F_{\rm pren} =\Flin \left ( \di\sum_{\bn \in A(F)^*} \Pren^{\bn} B_{\bn} \right ) ,
\end{equation}
where $A(F)^*$ is the set of sequences $\bn =(n_1 ,\dots ,n_r)$, $n_i \in A(F)$, $r\geq 0$, $\Pren^{\bn} \in \C$ satisfying
\begin{equation}
\Pren^{\bn} =0\ \ \ \mbox{\rm if}\ \ \parallel \bn \parallel \not= 0 ,
\end{equation}
with $\parallel \bn \parallel =n_1 +\dots +n_r$ for all $\bn \in A(F)^*$, $\bn=n_1 \dots n_r$, and $B_{\bn} =B_{n_1} \dots B_{n_r}$ with the usual composition of differential operators.
\end{defi}

These forms are calculable using the formalism of {\it moulds} developed by J. Ecalle since 1970.

\section{Moulds and prenormalization}

\subsection{Reminder about moulds}

We provide a self-contained introduction to the formalism of moulds and we refer to the articles of J. Ecalle or to the lectures
(\cite{cr1},\cite{cr2}) for more details.

\subsubsection{Moulds and non-commutative formal power series}

We denote by $A$ an alphabet, finite or not. A letter of $A$ is denoted by $a$. Let $A^*$ denotes the set of {\it words}
constructed on $A$, {\it i.e.} the sequences $a_1 \dots a_r$, $r\geq 0$, with $a_i \in A$, with the convention that for $r=0$
we have the {\it empty-word} denoted by $\emptyset$. We denote with bold letter $\ba$ a word of $A^*$. We have a natural action on $A^*$ provided by the usual
{\it concatenation} of two words $\ba$, $\bb \in A^*$, which glues the words $\ba$ to $\bb$, {\it i.e.} $\ba \bb$.

\begin{defi}
Let $\K$ be a ring (or a field) and $A$ a given alphabet. A $\K$-valued mould on $A$ is a map from $A^*$ to $\K$, denoted
by $\Mb$.
\end{defi}

The evaluation of $\Mb$ on a word $\ba \in A^*$ is denoted by $\M^{\ba}$\\

As an example, we define a $\C$-valued mould on $A(F)$ by
\begin{equation}
\left .
\begin{array}{llll}
\Prenb : & A(F)^* & \longrightarrow & \C \\
 & \bn & \longmapsto & \Pren^{\bn} .
\end{array}
\right .
\end{equation}

The mould $\Prenb$ is obtained collecting the coefficients of a formal power serie
$\di\sum_{\bn \in A(F)^*} \Pren^{\bn} B_{\bn}$. There exist a one-to-one correspondence between moulds and formal
power series.\\

For $r\geq 0$, we denote by $A^*_r$ the set of words of length $r$, with the convention that
$A^*_0 =\{ \emptyset \}$. We denote by $\K \langle A\rangle$ the set of finite $\K$-linear combinations of elements of $A^*$, {\it i.e.}
{\it non-commutative} polynomials on $A$ with coefficients in $K$, and by $\K_r \langle A\rangle$ the
set of $\K$-linear combination of elements of $A^*_r$, {\it i.e.} the set of non-commutative homogeneous polynomials of
degree $r$. We have a natural {\it graduation} on $\K \langle A\rangle$ by the length of words:
\begin{equation}
\K \langle A\rangle =\di\bigoplus_{r=0}^{\infty} \K_r \langle A\rangle .
\end{equation}
The completion of $\K \langle A\rangle$ with respect to the graduation by length denoted by $\K \langle\langle
A\rangle\rangle$ is the set of formal power series with coefficients in $\K$. An element of $\K \langle\langle
A\rangle\rangle$ is denoted by
\begin{equation}
\di\sum_{\ba \in A^*} M^{\ba} \ba ,\ \ M^{\ba} \in \K ,
\end{equation}
where this sum must be understood as
\begin{equation}
\di\sum_{r\geq 0} \left ( \di\sum_{\ba\in A^*_r } M^{\ba} \ba \right ) .
\end{equation}
Let $\Mb$ be a $\K$-valued mould on $A$, its generating serie denoted by $\Phi_M$ belongs to $\K \langle\langle A\rangle\rangle$ and is defined by
\begin{equation}
\Phi_M =\di\sum_{\ba \in A^*} M^{\ba} \ba ,
\end{equation}
or in a condensed way as $\di\sum_{\bullet} M^{\bullet} \bullet $. This correspondence provide a {\it one-to-one mapping} from the set of $\K$-valued moulds
on $A$ denoted by ${\cal M}_{\K} (A)$ and $\K \langle\langle A\rangle\rangle$.

\subsubsection{Moulds algebra}

The set of moulds ${\cal M}_{\K} (A)$ inherits a {\it structure of algebra} from $\K \langle\langle A\rangle\rangle$. The sum and
product of two moulds $\Mb$ and $\Nb$ is denoted by $\Mb +\Nb$ and $\Mb \cdot \Nb$ respectively and defined by
\begin{equation}
\left .
\begin{array}{lll}
(\Mb +\Nb )^{\ba} & = & M^{\ba} +N^{\ba} ,\\
(\Mb \cdot \Nb )^{\ba} & = & \di\sum_{\ba^1 \ba^2 =\ba} M^{\ba^1} N^{\ba^2} ,
\end{array}
\right .
\end{equation}
for all $\ba \in A^*$ where the sum corresponds to all the partition of $\ba$ as a concatenation of two words $\ba^1$ and $\ba^2$ of $A^*$.\\

The neutral element for the mould product is denoted by $\1b$ and defined by
\begin{equation}
\1b =\left \{
\begin{array}{ll}
1& \ \ \mbox{\rm if}\ \bullet =\emptyset ,\\
0& \ \ \mbox{\rm otherwise},
\end{array}
\right .
\end{equation}

Let $\Mb$ be a mould. We denote by $\Mib$ the inverse of $\Mb$ for the mould product when it exists, {\it i.e.} the solution of the mould
equation:
\begin{equation}
\Mb \cdot \Mib =\Mib \cdot \Mb =\1b .
\end{equation}

\subsubsection{Composition of moulds}

Assuming that $A$ possesses a {\it semi-group} structure, we can define a non-commutative version of the classical operation of {\it substitution} of
formal power series.\\

We denote by $\star$ an internal law on $A$, such that $(A ,\star )$ is a semi-group. We denote by
$\parallel \, \parallel_{\star}$ the mapping from $A^*$ to $A$ defined by
\begin{equation}
\left .
\begin{array}{llll}
\parallel \, \parallel_{\star} :\ \ & A^* & \longrightarrow & A ,\\
 & \ba =a_1 \dots a_r & \longmapsto & a_1 \star\dots \star a_r .
\end{array}
\right .
\end{equation}
The $\star$ will be omitted when clear from the context.\\

The set $\K\langle\langle A\rangle\rangle$ is graded by $\parallel \, \parallel_{\star}$. A {\it homogeneous component} of degree $a\in A$ of a
non-commutative serie $\Phi_{M} =\di\sum_{\ba \in A^*} M^{\ba} \ba$ is the quantity
\begin{equation}
\Phi_M^a =\di\sum_{\ba \in A^*,\ \parallel \ba \parallel_{\star} =a} M^{\ba} \ba .
\end{equation}
We have by definition
\begin{equation}
\Phi_M =\di\sum_{a\in A} \Phi_M^a .
\end{equation}

\begin{defi}[Composition]
Let $(A ,\star )$ be a semi-group structure. Let $\Mb$ and $\Nb$ be two moulds on ${\cal M}_{\K} (A)$ and $\Phi_M$, $\Phi_N$ their associated
generating series. The substitution of $\Phi_N$ in $\Phi_M$, denoted by $\Phi_M \circ \Phi_N$ is defined by
\begin{equation}
\label{substi2}
\Phi_M \circ \Phi_N =\di\sum_{\ba \in A^*} M^{\ba} \Phi_N^{\ba} ,
\end{equation}
where $\Phi_N^{\ba}$ is given by $\Phi_N^{a_1} \dots \Phi_N^{a_r}$ for $\ba =a_1 \dots a_r$.\\

We denote by $\Mb \circ \Nb$ the mould of ${\cal M}_{\K} (A)$ such that
\begin{equation}
\label{substi}
\Phi_M \circ \Phi_N =\di\sum_{\ba \in A^*} (\Mb \circ \Nb )^{\ba} \ba .
\end{equation}
\end{defi}

Equation (\ref{substi}) define a natural operation on moulds denoted $\circ$ and called {\it composition}. Using $\parallel \, \parallel_{\star}$ we can
give a closed formula for the composition of two moulds.

\begin{lem}
Let $(A,\star )$ be a semi-group and $\Mb$, $\Nb$ be two moulds of ${\cal M}_{\K} (A)$. We have for all $\ba \in A^*$,
\begin{equation}
\label{formulacompo}
(\Mb \circ \Nb )^{\ba} =\di\sum_{k=1}^{l(\ba )} \sum_{\ba^1 \dots \ba^k \stackrel{=}{*} \ba } M^{\parallel \ba^1 \parallel_{\star}
\dots \parallel \ba^k \parallel_{\star}} N^{\ba^1} \dots N^{\ba^k} ,
\end{equation}
where $\ba^1 \dots \ba^k \stackrel{=}{*} \ba$ denotes all the partitions of $\ba$ such that $\ba^i \not= \emptyset$, $i=1,\dots ,k$.
\end{lem}

\begin{proof}
Equation (\ref{substi2}) is equivalent to
\begin{equation}
\label{sommecompo}
\Phi_M \circ \Phi_N =\di\sum_{r\geq 0} \sum_{\bb =b_1 \dots b_r \in A^*_r} M^{b_1 \dots b_r} \left ( \di\sum_{\ba^1 \in A^*,\ \parallel \ba^1 \parallel_{\star} =b_1} N^{\ba^1} \ba^1 \right )
\dots \left ( \di\sum_{\ba^r \in A^*,\ \parallel \ba^r \parallel_{\star} =b_r} N^{\ba^r} \ba^r \right ) .
\end{equation}
Let $\ba \in A^*$ be a given word of $A^*$. Each partition of $\ba$ of the form $\ba =\ba^1 \dots \ba^k$, $k=1,\dots ,l(\ba )$, occurs in the sum
(\ref{sommecompo}) with a coefficient given by
\begin{equation}
M^{b_1 \dots b_r} N^{\ba^1} \dots N^{\ba^k} ,
\end{equation}
where $b_i =\parallel \ba^i \parallel_{\star}$. Collecting all these coefficients, we obtain the formula (\ref{formulacompo}) for the coefficient of $\ba$
in $\Phi_M \circ \Phi_N$.
\end{proof}

The neutral element for the mould composition is denoted by $\Ib$ and defined by
\begin{equation}
\Ib =\left \{
\begin{array}{ll}
1& \ \ \mbox{\rm if}\ l(\bullet) =1 ,\\
0& \ \ \mbox{\rm otherwise},
\end{array}
\right .
\end{equation}
where $l(\bullet )$ denotes the length of a word of $A^*$.

\subsubsection{Exponential and logarithm of moulds}

We denote by $(\K \langle\langle A\rangle\rangle )_*$ the set of formal power series without a constant term. We define the {\it exponential} of an element $x\in (\K \langle\langle A\rangle\rangle )_*$
, denoted by $\exp (x)$ using the classical formula
\begin{equation}
\exp (x)=\di\sum_{n\geq 0} \di {x^n \over n!} .
\end{equation}
The {\it logarithm} of an element $1+x\in 1+(\K \langle\langle A\rangle\rangle )_*$ is denoted by $\log (1+x)$ and defined by
\begin{equation}
\log (1+x )=\di\sum_{n\geq 0} (-1)^{n+1} \di {x^n \over n!} .
\end{equation}
These two applications have their natural counterpart in ${\cal M}_{\K} (A)$.

\begin{defi}
Let $\Mb$ be a mould of ${\cal M}_{\K} (A)$ and $\Phi_M$ the associated generating serie. Assume that $\exp (\Phi_M )$ is defined. We denote by
$\Exp \Mb$ the mould satisfying the equality
\begin{equation}
\exp \left( \di\sum_{\bullet} \Mb \bullet \right ) =\di\sum_{\bullet} \Exp \Mb \bullet .
\end{equation}
\end{defi}

Simple computations lead to the following direct definition of $\Exp$ on moulds:
\begin{equation}
\Exp \Mb =\di\sum_{n\geq 0} \di {\left [ \Mb \right ]_{(\times n)} \over n!} ,
\end{equation}
where $\left [ \Mb \right ]_{(\times n)}$, $n\in \N$, stands for
\begin{equation}
\left [ \Mb \right ]_{(\times n)} =\underbrace{\Mb \cdots \Mb}_{n\ \mbox{\rm times}} .
\end{equation}

The same procedure can be applied to define the logarithm of a mould.

\begin{defi}
Let $\Mb$ be a mould of ${\cal M}_{\K} (A)$ and $\Phi_M$ the associated generating serie. Assume that $\log (1+\Phi_M )$ is defined. We denote by
$\Log \Mb$ the mould satisfying the equality
\begin{equation}
\log \left( 1+ \di\sum_{\bullet} \Mb \bullet \right ) =\di\sum_{\bullet} \Log \Mb \bullet .
\end{equation}
\end{defi}

A direct definition of $\Log$ is then given by
\begin{equation}
\Log \Mb =\di\sum_{n\geq 0} (-1)^{n+1} \di {\left [ \Mb \right ]_{(\times n)} \over n!} .
\end{equation}

As $\exp$ and $\log$ satisfy $\exp \circ \log =\log \circ \exp =1$, we have
\begin{equation}
\Exp \left ( \Log \Mb \right ) =\Log \left ( \Exp \Mb \right ) =\1b .
\end{equation}

\subsubsection{A technical lemma}

In this section, we derive simple results for the exponential and logarithm of moulds with non-zero components only on words of length $1$.\\

\begin{lem}
\label{simple}
Let us denote by $\Zib$ a mould of ${\cal M}_{\K} (A)$ such that $\Zib =0$ for all $\bullet$ of length different from $1$. For all $\ba \in A^*$,
$r\geq 1$, we have
\begin{eqnarray}
\left [ \Zib \right ]_{\times r}^{\ba} =\left \{
\begin{array}{l}
\Zi^{a_1} \dots \Zi^{a_2} \ \mbox{\rm if}\ \ l(\ba )=r,\ \ba =a_1 \dots a_r, \\
0\ \ \mbox{\rm otherwise.}
\end{array}
\right .
\label{prodz}\\
\left [\Exp \Zib \right ]^{\ba} =\un^{\ba} +\di {1\over l(\ba )!} \left [\Zib \right ]_{(\times l(\ba ))}^{\ba} , \label{expoz}\\
\left [\Log \Zib \right ]^{\ba} = \di {(-1)^{l(\ba )+1} \over l(\ba )!} \left [ \Zib \right ] _{(\times l(\ba ))}^{\ba} \label{logz}.
\end{eqnarray}
\end{lem}

\begin{proof}
We first remark that equations (\ref{expoz}) and (\ref{logz}) easily follow from equation (\ref{prodz}).\\

The proof of equation (\ref{prodz}) is done by induction on $r$. Formula (\ref{prodz}) is trivially true for $r=1$. Assume that formula (\ref{prodz})
is true for $r\geq 1$. By definition, we have
\begin{equation}
\left [ \Zib \right ]_{(\times r+1)} =\Zib \times \left [ \Zib \right ]_{(\times r)} .
\end{equation}
Let $\ba =a \bb$, then by assumption on $\Zib$ we obtain
\begin{equation}
\left [ \Zib \right ]_{(\times r+1)}^{a\bb} =\Zi^a \left [ \Zib \right ]_{(\times r)}^{\bb} .
\end{equation}
As the mould $\left [ \Zib \right ]_{(\times r)}$ is non-trivial only on words of length $r$, we deduce that the mould
$\left [ \Zib \right ]_{(\times r+1)}$ is non-trivial only on words of length $r+1$.

Moreover, using the fact that $\left [ \Zib \right ]_{(\times r)}^{a_1 \dots a_r} =\Zi^{a_1} \dots \Zi^{a_r}$ for all $a_i \in A$, we also
deduce that $\left [ \Zib \right ]_{(\times r+1)}^{a_1 \dots a_{r+1}} =\Zi^{a_1} \dots \Zi^{a_{r+1}}$. This concludes the proof.
\end{proof}

\subsection{Prenormalization}

Let $F$ be a diffeo in prepared form given by
$$F=\Flin \left ( \mbox{\rm Id} +\di\sum_{n\in A(F)} B_n \right ) .$$
Let $\Phi_{\Theta}$ be an automorphism of $\C \{ x\}$ of the form
\begin{equation}
\Phi_{\Theta} =\di\sum_{\bn \in A(F)^*} \Theta^{\bn} B_{\bn} ,
\end{equation}
where $\Theta^{\bn} \in \C$ for all $\bn\in A(F)^*$, {\it i.e.} $\Phi_{\Theta} \in \C \langle\langle \Bbb \rangle\rangle$,
where $\Bbb =\{ B_n \}_{n\in A(F)}$ and $\Thetab \in {\cal M}_{\C} (A(F))$.\\

Using the moulds $\1b$ and $\Ib$ we write $\mbox{\rm Id} +\di\sum_{n\in A(F)} B_n$ as an element of $\C \langle\langle \Bbb \rangle\rangle$:
\begin{equation}
\mbox{\rm Id} +\di\sum_{n\in A(F)} B_n =\di\sum_{\bullet} (\1b +\Ib ) \Bb .
\end{equation}
We assume that $F$ is conjugated to an automorphism $F_{\rm conj}$ via $\Phi_{\Theta}$. Equation (\ref{conjugacy}) is then given by
\begin{equation}
\label{conjugacy2}
F_{\rm conj} = \Phi_{\Theta} \cdot F \cdot \Phi_{\Theta}^{-1} .
\end{equation}
The automorphism $F_{\rm conj}$ can be written as
\begin{equation}
F_{\rm conj} =\Flin \left ( \di\sum_{\bullet} \Conjb \Bb \right ).
\end{equation}
Equation (\ref{conjugacy2}) is then equivalent to
\begin{equation}
\label{conjugacy3}
\Flin \left  ( \di\sumb \Conjb \Bb \right ) =\left ( \di\sumb \Thetab \Bb \right ) \Flin \left ( \di\sumb (\1b +\Ib )\Bb \right )
\left ( \di\sumb \Thetaib \Bb \right ) ,
\end{equation}
where $\Thetaib$ is such that $\Thetaib \cdot \Thetab =\Thetab \cdot \Thetaib =\1b$, {\it i.e.} $\Phi_{\Theta}^{-1} =\di\sumb \Thetaib \Bb$.\\

In order to explicit $\Conjb$ we need to understand the action of a formal power serie of $\C \langle \langle \Bbb \rangle\rangle$ on $\Flin$.
We have the following fundamental lemma:

\begin{lem}
\label{actionlin}
Let $\Mb \in {\cal M}_{\C} (A(F))$ . We have
\begin{equation}
\left ( \di\sumb \Mb \Bb \right ) \Flin =\Flin \left ( \di\sumb \edelta \left (\Mb \right )^\bullet \Bb \right ) ,
\end{equation}
where $\edelta$ is a map from ${\cal M}_{\C} (A(F))$ to ${\cal M}_{\C} (A(F))$ defined by
\begin{equation}
\edelta \left (\Mb \right )^{\bn} =e^{-\lambda .\parallel \bn \parallel} M^{\bn} \ \mbox{\rm for all}\ \bn\in A(F)^* .
\end{equation}
\end{lem}

\begin{proof}
Let $B_{\bn} =B_{n_1 \dots n_r}$ such that $B_{n_i} (x^m )=\beta^{n_i}_m x^{m+n_i}$, $\beta^{n_i}_m \in \C$, $i=1,\dots ,r$, for all $m\in \N^{\nu}$.
We have
\begin{equation}
B_{\bn} (x^m )=\beta^{n_1}_{m+n_r +\dots +n_2} \beta^{n_2}_{m+n_r +\dots +n_3} \dots \beta^{n_r}_m x^{m+n_1 +\dots +n_r} .
\end{equation}
As $\Flin (x^m )=\di e^{\lambda .m} x^m$ we obtain
\begin{equation}
\left .
\begin{array}{lll}
B_{\bn} \left ( \Flin (x^m ) \right ) & = & \di e^{\lambda .m} B_{\bn} (x^m ) ,\\
 & = & \di e^{-\lambda .(n_1 +\dots +n_r )} \di e^{\lambda . (m+n_1 +\dots +n_r )} B_{\bn} (x^m ) ,\\
  & = & \di e^{-\lambda .(n_1 +\dots +n_r )} \Flin \left ( B_{\bn} (x^m ) \right ) ,\\
  & = & \di \Flin \left ( e^{-\lambda .(n_1 +\dots +n_r )} B_{\bn} (x^m ) \right ) .
\end{array}
\right .
\end{equation}
This concludes the proof.
\end{proof}

Next lemma gives an explicit formula to compute the mould $\Conjb$ assuming that the mould $\Thetab$ is known.

\begin{lem}
\label{important}
Equation (\ref{conjugacy3}) is equivalent to the mould equation
\begin{equation}
\label{formuleimportante}
\Conjb =\di\edelta \left ( \Thetab \right ) \cdot (\1b +\Ib ) \cdot \Thetaib .
\end{equation}
\end{lem}

\begin{proof}
Using lemma \ref{actionlin}, we have
\begin{equation}
\left .
\begin{array}{lll}
\Flin \left  ( \di\sumb \Conjb \Bb \right ) & = & \left ( \di\sumb \Thetab \Bb \right ) \Flin \left ( \di\sumb (\1b +\Ib )\Bb \right )
\left ( \di\sumb \Thetaib \Bb \right ) ,\\
 & = & \Flin \left ( \di\sumb \edelta \left ( \Thetab \right ) \Bb \right ) \left ( \di\sumb (\1b +\Ib )\Bb \right )
\left ( \di\sumb \Thetaib \Bb \right ) ,\\
 & = & \Flin \left ( \di\sumb \left ( \edelta \left ( \Thetab \right ) \cdot (\1b +\Ib ) \cdot \Thetaib \right )  \Bb \right ) .
\end{array}
\right .
\end{equation}
This concludes the proof.
\end{proof}

As a consequence, choosing carefully the normalizator $\Phi_{\Theta}$, we can obtain an inductive expression for the mould of normalization $\Conjb$.\\

We will give explicit formulae for $\Conjb$ using specific moulds for $\Thetab$ in the next section.

\subsection{Universality of moulds and prenormalization}

Lemma \ref{important} gives an important feature of the mould formalism in the context of continuous prenormalization. Formula (\ref{formuleimportante}) is
valid whatever is the underlying alphabet $A(F)$. We then obtain a {\it universal} object underlying the prenormalization problem which is studied.\\

For example, in the context of {\it linearization}, {\it i.e.} $F_{\rm conj} =\Flin$, the universal mould of linearization which defined the linearizing
change of variables is given as follow (see \cite{cr1} Chap. III for more details):

\begin{thm}
Let $\mbox{\rm\bf L} =\{ L_r \}_{r\geq 1}$, $r\in \N$, be the set of $\C$-valued functions $L_r :\C^r \rightarrow \C$ defined by
\begin{equation}
L_r (x_1 ,\dots ,x_r )=\left [ \left ( \di e^{-(x_1 +\dots +x_r )} -1 \right ) \left ( \di e^{-(x_2 +\dots +x_r )} -1 \right ) \dots \left ( \di e^{-x_r } -1 \right )
\right ] ^{-1} ,
\end{equation}
for all $(x_1 ,\dots ,x_r )\in \C^r \setminus S_r$ where the singular set $S_r$ is given by
\begin{equation}
S_r =\{ x_r =0\} \bigcup \{ x_r +x_{r-1} =0 \} \bigcup \dots \bigcup \{ x_1 +\dots +x_r =0 \} .
\end{equation}
If $F$ possesses a non-resonant linear part $\lambda$, the mould of formal linearization is given for all $\bn \in A(F)^*$, $\bn =n_1 ,\dots ,n_r$, by
\begin{equation}
\label{linearization}
\Theta^{n_1 \dots n_r} =L_r (\omega_1 ,\dots ,\omega_r ) ,
\end{equation}
where $\omega_i =n_i .\lambda$ for $i=1,\dots ,r$.
\end{thm}

This result can not be obtained using other existing formalisms. Of course, anybody knows that an expression like (\ref{linearization}) is the important
quantity entering the linearization problem. However, the previous result associates universal coefficients from which one can compute the desired
linearization map for a given particular diffeo $F$ by posing
$$
\Phi_{\Theta} =\sum_{\bn \in A(F)^*} \Theta^{\bn} B_{\bn} .
$$

\section{The Trimmed form}

In this section, we give complete proofs of results concerning the {\it Trimmed form} defined by J. Ecalle and B. Vallet announced in \cite{ev1} with a
sketch of proof and without explicit computations.

\subsection{Cancelling non-resonant terms}
\label{cancel}

In this section, we give a mould approach to the classical problem of cancellation of non-resonant terms.

\subsubsection{Around the Baker-Campbell-Hausdorff formula}

Let $F$ be a diffeo in prepared form given by (\ref{homocomponents}). The operator $\mbox{\rm Id}+\di\sum_{n\in A(F)} B_n$ is an automorphism of
$\C \{ x\}$ which can be viewed as the exponential of a vector field, {\it i.e.}
\begin{equation}
\label{vectordecompo}
\mbox{\rm Id}+\di\sum_{n\in A(F)} B_n =\exp \left ( \di\sum_{m\in {\cal A}(F)} D_m \right ) ,
\end{equation}
where $D_n$ is a homogeneous differential operator of degree $m$ and order $1$, {\it i.e.} a derivation on $\C \{ x\}$, $m=(m_1 ,\dots ,m_{\nu} )\in \Z^{\nu}$, with all $m_i \in \N$,
$i=1,\dots ,\nu$ except at most one which can be $-1$, and ${\cal A} (F)$ the set of degrees coming in the decomposition.\\

We look for an automorphism given by the exponential of a vector field $\mbox{\rm\bf V}$ given by
\begin{equation}
\mbox{\rm\bf V} =\di\sum_{\bn \in A(F)^*} \dem^{\bn} B_{\bn} ,
\end{equation}
or equivalently given on the alphabet ${\cal A}(F)^*$ by
\begin{equation}
\mbox{\rm\bf V} =\di\sum_{\bbm \in {\cal A}(F)^*} \Dem^{\bbm} D_{\bbm} ,
\end{equation}
where $V_{\bbm =m_1 \dots m_r} =V_{m_1} V_{m_2} \dots V_{m_r}$ with the usual composition of differential operators.\\

The action of $\exp \mbox{\rm\bf V}$ on $F$ is given by
\begin{equation}
\label{actionV}
\exp \mbox{\rm\bf V} \cdot F \cdot \exp (-\mbox{\rm\bf V} )
\end{equation}
Equation (\ref{actionV}) can be analyzed using the moulds expression of $\mbox{\rm\bf V}$ and $F$ with respect to the alphabet ${\cal A} (F)$. We
have the following lemma:

\begin{lem}
\label{lemmeactionV}
Equation (\ref{actionV}) is equal to
\begin{equation}
\exp \mbox{\rm\bf V} \cdot F \cdot \exp (-\mbox{\rm\bf V} )
= \Flin \exp \left ( \tilde{\mbox{\rm\bf V}} +\mbox{\rm\bf D} -\mbox{\rm\bf V}+\dots \right ) ,
\end{equation}
where the $\dots$ stands for a formal power serie beginning with words of length at least $2$, and $\mbox{\rm\bf D}$ and $\tilde{\mbox{\rm\bf V}}$ are
vector fields defined by $\mbox{\rm\bf D} =\di\sum_{m\in {\cal A} (F)} D_m$ and
\begin{equation}
\tilde{\mbox{\rm\bf V}} = \di\sum_{\bbm\in {\cal A} (F)^*} \di e^{-\lambda .\parallel\bbm\parallel} \Dem^{\bbm} D_{\bbm} ,
\end{equation}
respectively.
\end{lem}

\begin{proof}
Using the Baker-Campbell-Hausdorff formula ($BCH_2$), we obtain
\begin{equation}
\left .
\begin{array}{lll}
\exp \mbox{\rm\bf D} \cdot \exp (-\mbox{\rm\bf V} ) & = & \exp \left ( \mbox{\rm\bf D} \star (-\mbox{\rm\bf V} ) \right ) ,\\
 & = & \exp \left(\mbox{\rm\bf D} -\mbox{\rm\bf V} +\mbox{\rm h.o.t.}\right) ,
\end{array}
\right .
\end{equation}
where $\mbox{\rm h.o.t.}$ stands for {\it higher order terms}.\\

Using lemma \ref{actionlin}, we have
\begin{equation}
\exp \mbox{\rm\bf V} \cdot \Flin =\Flin \cdot \exp \tilde{\mbox{\rm\bf V}},
\end{equation}
where $\tilde{\mbox{\rm\bf V}}$ is given by
\begin{equation}
\tilde{\mbox{\rm\bf V}} = \di\sum_{\bbm\in {\cal A} (F)^*} \di e^{-\lambda .\parallel\bbm\parallel} \Dem^{\bbm} D_{\bbm} .
\end{equation}
As a consequence, applying again ($BCH_2$) we obtain
\begin{equation}
\left .
\begin{array}{lll}
\exp \tilde{\mbox{\rm\bf V}} \cdot \exp (\mbox{\rm\bf D} \star (-\mbox{\rm\bf V}) ) & = & \exp \left ( \tilde{\mbox{\rm\bf V}} \star
(\mbox{\rm\bf D} \star (-\mbox{\rm\bf V} ) ) \right ) ,\\
 & = & \exp \left ( \tilde{\mbox{\rm\bf V}} +\mbox{\rm\bf D} -\mbox{\rm\bf V} +\dots \right ) ,
\end{array}
\right .
\end{equation}
where the $\dots$ stand for a formal power serie beginning with words of length at least $2$. This concludes the proof.
\end{proof}

\subsubsection{The simplified form and the moulds $\demb$ and $\Demb$}
\label{sectionsimplified}

The main consequence of lemma \ref{lemmeactionV} is that we can cancel the non-resonant terms of $\mbox{\rm\bf D}$ using a simple vector field
$\mbox{\rm\bf V}$.

\begin{thm}[Simplified form]
\label{tsimplifiedform}
Let $\mbox{\rm\bf V}$ be the vector field defined by the mould
\begin{equation}
\Demb =\left \{
\begin{array}{l}
\di {I^{\bbm} \over 1-\di e^{\mid \bbm \mid .\lambda} } \ \ \mbox{\rm for}\ \ \bbm \in {\cal A} (F)^* \setminus {\cal R} (F) ,\\
0\ \ \mbox{\rm otherwise},
\end{array}
\right .
\end{equation}
where ${\cal R}(F)$ is the set of resonant words of ${\cal A} (F)^*$, {\it i.e.} $\bbm \in {\cal R} (F)$ if and only if $\bbm .\lambda =0$. We
denote by $\demb$ the associated mould on ${\cal M}_{\C} (A(F))$, {\it i.e.}
\begin{equation}
\mbox{\rm\bf V} =\di\sumb \Demb \Db =\di\sumb \demb \Bb .
\end{equation}
We call simplified form of $F$ and we denote by $\FSem$ the automorphism obtained from $F$ under the action of $\exp \mbox{\rm\bf V}$. We have
\begin{equation}
\left .
\begin{array}{lll}
\FSem & = & \Flin \left ( \di\sum_{\bbm \in {\cal A} (F)^*} \Sem^{\bbm} D_{\bbm} \right ) ,\\
      & = & \Flin \left ( \di\sum_{\bn \in A (F)^*} \sem^{\bn} B_{\bn} \right )
\end{array}
\right .
\end{equation}
with the mould $\Semb$ given by
\begin{equation}
\Semb =\edelta \left ( \Exp (\Demb )\right ) \cdot \Exp (\Ib )\cdot \Exp (-\Demb ) ,
\end{equation}
and the mould $\semb$ given by
\begin{equation}
\semb =\edelta \left ( \Exp (\demb )\right ) \cdot (\1b +\Ib)\cdot \Exp (-\demb ) .
\end{equation}
\end{thm}

\begin{proof}
We have $\FSem =\exp \mbox{\rm\bf V} \cdot F \cdot \exp (-\mbox{\rm\bf V} )$ with $\mbox{\rm\bf V}=\di\sum_{\bn \in A (F)^*} \dem^{\bn} B_{\bn}$. As
a consequence, we have $\exp \mbox{\rm\bf V} =\di\sum_{\bn \in A(F)^*} \left ( \Exp \, \demb \right )^{\bn} B_{\bn}$ and the
formula for $\semb$ follows from lemma \ref{important} using $\Thetab =\Exp (\demb )$.\\

For $\Semb$, we first use lemma \ref{actionlin} to obtain
\begin{equation}
\exp \mbox{\rm\bf V} \Flin =\Flin \di \left ( \di\sum_{\bbm \in {\cal A} (F)^*} \left [ \edelta \left ( \Exp \Demb \right ) \right ] ^{\bn}
D_{\bbm} \right ) .
\end{equation}
As a consequence, the conjugacy equation is equivalent to
\begin{equation}
\left .
\begin{array}{lll}
\FSem & = & \exp \mbox{\rm\bf V} \cdot F \cdot \exp (-\mbox{\rm\bf V} ) ,\\
   & = & \Flin \left ( \di\sumb \edelta \left ( \Exp \Demb \right ) \Db \right ) \left ( \di\sumb \Exp \Ib \Db \right )
   \left ( \di\sumb \Exp ( -\Demb ) \Db \right ) ,\\
   & = & \Flin \left ( \di\sumb \left [ \edelta \left ( \Exp \Demb \right ) \cdot \Exp \Ib \cdot \Exp (-\Demb ) \right ]^\bullet \Db \right ) .
\end{array}
\right .
\end{equation}
This concludes the proof.
\end{proof}

The mould $\Semb$ can be compute explicitly. We first introduce some convenient notations:\\

Let $\bbm =m_1 \dots m_r$ be a word of length $r$, $r\geq 1$. We denote by $\bbm^{\leq i}$ and $\bbm^{>i}$ the word
\begin{equation}
\bbm^{\leq i} =m_1 \dots m_i ,\ \ \ \bbm^{>i} =m_{i+1} \dots m_r .
\end{equation}
Moreover we denote by $d(\bbm)$ the index of the last $m_i$ in $\bbm = m_1 \dots m_r$ such that $\lambda. m_i = 0$, and we denote by $q(\bbm)$ the first index just before of the first zero $\omega_j = \lambda. m_j$.

\begin{thm}
\label{explicitSem}
For all $\bbm \in {\cal A} (F)^*$, we have
\begin{equation}
\small
\left .
\begin{array}{l}
\Sem^{\bbm}  =  \di {(-1)^{l(\bbm)}\over l(\bbm )!}\left [\Demb \right ]^{\bbm}_{(\times l(\bbm ))} +\di {1\over l(\bbm )!} +\di\sum_{j=d(\bbm )+1}^{l(\bbm )+1}
\di { (-1)^{l(\bbm^{>j})}\left [\Demb \right ]^{\bbm^{\geq j}}_{(\times l(\bbm^{\geq j} ))} \over l(\bbm^{<j} )! l(\bbm^{\geq j} )! } +\di e^{-\lambda. \parallel \bbm \parallel} \, \mbox{\rm 1}^{\bbm} \\
  +\di\sum_{i=1}^{q(\bbm ) \wedge (l(\bbm)-1)} \di {e^{-\lambda. \parallel \bbm ^{\leq i}\parallel} \over l(\bbm^{\leq i} )!} \left [\Demb \right ]
^{\bbm^{\leq i}}_{(\times l(\bbm^{\leq i} ))} \times \\
 \left (
 \di  {(-1)^{l(\bbm)}\over l(\bbm )!} \left [ \Demb \right ]^{\bbm^{>i}}_{(\times l(\bbm^{>i} ))} +\di {1\over l(\bbm^{>i} )!} +
 \di\sum_{j=d(\bbm^{>i} )+1}^{l(\bbm^{>i} )+1}
 \di {(-1)^{l(\bbm^{>j})^{\geq j}} \left [\Demb \right ] ^{(\bbm^{>i} )^{\geq j}}_{(\times l((\bbm^{>i} )^{\geq j} ))} \over
 l((\bbm^{>i} )^{<j} )! l((\bbm^{>i})^{\geq j} )!} \right ) .
\end{array}
\right .
\end{equation}
\end{thm}

The proof is done in appendix \ref{proofSem}.

\subsection{The Trimmed form}

The Trimmed form is constructed by induction applying successively the previous simplification scheme to remove non-resonant terms of higher and
higher degrees. The mould formalism allows us to explicit some particular moulds underlying this construction as well as algorithmic and
explicit formulae for some of them.

\subsubsection{The Trimmed form up to order r}

We can use the simplification procedure previously defined inductively in order to cancel non-resonant terms of higher and higher degrees.

\begin{defi}[Trimmed form up to order $r$]
Let $r\in \N$, the Trimmed form up to order $r$ is defined as $\FSem^r$ obtained from $F$ after $r$ successive simplifications, {\it i.e.}
\begin{equation}
F=\FSem^0 \stackrel{\AutSimp^1}{\rightarrow} \FSem^1 \stackrel{\AutSimp^2}{\rightarrow} \dots \stackrel{\AutSimp^r}{\rightarrow} \FSem^r ,
\end{equation}
where $\AutSimp^i$ is the automorphism of simplification defined by
\begin{equation}
\AutSimp^i =\exp ( \mbox{\rm\bf V}_i ) ,
\end{equation}
with $\mbox{\rm\bf V}_i$ the vector fields associated to the mould $\Demb$ on the alphabet ${\cal A} (\FSem^{i-1} )$ associated to $\FSem^{i-1}$.
\end{defi}

Using theorem \ref{tsimplifiedform}, we deduce the following useful result:

\begin{thm}
For all $r\in \N$, the Trimmed form up to order $r$ denoted $\FSem^r$ possesses a mould expansion, {\it i.e.} there exists moulds denoted by
${}_r \Semb \in {\cal M}_{\C} ({\cal A}(F))$ and ${}_r \semb \in {\cal M}_{\C} (A(F))$ such that
\begin{equation}
\Fsem^r =\Flin \left ( \di\sumb {}_r \Semb \Db \right ) =\Flin \left ( \di \sumb {}_r \semb \Bb \right ) .
\end{equation}
\end{thm}

Despite its moulds expansion, the Trimmed form up to order $r$ is {\it not} a prenormal form as it remains non-resonant terms for sequences of
length $l\geq r+1$.

\subsubsection{The moulds ${}_r \semb$ and ${}_r \Semb$}

The mould ${}_r \semb$ has a simple expression in function of $\semb$.

\begin{lem}
For all $r\in \N$, we have
\begin{equation}
{}_r \semb =\underbrace{\semb \circ \dots \circ \semb}_{r\ \mbox{\rm times}} .
\end{equation}
\end{lem}

\begin{proof}
The simplification procedure can be written as follows:
\begin{equation}
\di\sumb \Ib \Bb \longmapsto \di\sumb \semb \Bb .
\end{equation}
Iterating this mapping we go from step $i$ to $i+1$
\begin{equation}
\di\sumb {}_i \semb \Bb =\di\sumb \Ib {} {}_{i+1} \Bb \longmapsto \di\sumb {}_{i+1} \semb  \Bb =\di\sumb \semb {}_{i+1} \Bb ,
\end{equation}
where $\di\sumb \Ib {}_{i+1} \Bb$ denotes the homogeneous decomposition constructed on $\Fsem^i$.\\

By definition of the composition for moulds we have
\begin{equation}
\di\sumb \semb {}_{i+1} \Bb =\di\sumb \left ( \semb \circ {}_i \semb \right ) \Bb ,
\end{equation}
from which we deduce the recursive relation
\begin{equation}
{}_{i+1} \semb =\semb \circ {}_i \semb .
\end{equation}
We conclude by induction on $i$.
\end{proof}

For the mould ${}_r \Semb$ we have a more complicated formula:

\begin{lem}
For all $r\in \N$, we have
\begin{equation}
\Log [{}_r \Semb ]=\underbrace{\Log (\Semb ) \circ \dots \circ \Log ( \Semb )}_{r\ \mbox{\rm times}} .
\end{equation}
\end{lem}

The fact that we must take the $\Log$ of $\Semb$ instead of $\Semb$ is related to the fact that the alphabet of derivation ${}_{i+1} \Db$ constructed
at step $i$ from $\Fsem^i$ is not related to $\sumb {}_i \Semb \Db$ but to its logarithm.

\begin{proof}
The simplification procedure can be written as follows:
\begin{equation}
\exp \left ( \di\sumb \Ib \Db \right ) \longmapsto \di\sumb \Semb \Db =\di\exp \left ( \di\sumb \Log (\Semb ) \Db \right ) .
\end{equation}
Iterating this mapping we go from step $i$ to $i+1$
\begin{equation}
\left .
\begin{array}{c}
\di\exp \left ( \di\sumb \Log [ {}_i \Semb ] \Db \right ) =\exp \left ( \di\sumb \Ib {} {}_{i+1} \Db \right ) \\
\downarrow \\
\exp \left ( \di\sumb \Log [{}_{i+1} \Semb ]  \Db \right ) =\exp \left ( \di\sumb \Log (\Semb ){}_{i+1} \Db \right ) ,
\end{array}
\right .
\end{equation}
where $\di\sumb \Ib {}_{i+1} \Db$ denotes the homogeneous decomposition constructed on $\di\sumb \Log [ {}_i \Semb ] \Db$.\\

By definition of the composition of moulds, we deduce that
\begin{equation}
\Log [{}_{i+1} \Semb ] =\Log (\Semb ) \circ \Log [ {}_i \Semb ] .
\end{equation}
We conclude the proof by induction on $i$.
\end{proof}

\subsubsection{The Trimmed form}

\begin{defi}
The Trimmed form of $F$ is the limit of the simplification procedure.
\end{defi}

\begin{thm}
The Trimmed form is a continuous prenormal form given by
\begin{equation}
\left .
\begin{array}{lll}
\FTrem & = & \Flin \left ( \di\sum_{\bbm \in {\cal A} (F)^*} \Trem^{\bbm} D_{\bbm} \right ),\\
       & = & \Flin \left ( \di\sum_{\bn \in A (F)^*} \trem^{\bn} B_{\bn} \right )
\end{array}
\right .
\end{equation}
with the moulds $\Tremb$ and $\tremb$ defined by
\begin{equation}
\begin{array}{c}
\Tremb -\1b = \limstat_{r\rightarrow \infty} \left  [ \Semb -\1b \right ] ^{(\circ r)} ,\\
\tremb -\1b =\limstat_{r\rightarrow \infty} \left  [ \semb -\1b \right ] ^{(\circ r)} ,
\end{array}
\end{equation}
where $\limstat$ is the stationary limit.
\end{thm}

The proof is a direct consequence of the simplification procedure.

\begin{rema}
Following (\cite{ev1} $\S$.7) we have divergence and resurgence of the simplification procedure. This is not the case when working directly
with the diffeomorphism instead of its associated automorphism of substitution. However, this problem can be avoided (see \cite{ev1} p.8).
\end{rema}

\subsubsection{The mould $\Tremb$}

We can compute the mould $\Tremb$ using a simple remark. By definition, we have the following identities
\begin{eqnarray}
\Tremb = \Semb \circ \Tremb ,\\
\Tremb = \Tremb \circ \Semb \label{etrem2}.
\end{eqnarray}
Using the first equation and the definition of composition for moulds we obtain for all $\bbm \in {\cal A} (F)^*$
\begin{equation}
\label{trem1}
\Trem^{\bbm} =\Sem^{\parallel \bbm \parallel} \Trem^{\bbm} +\mbox{\rm s.l} ,
\end{equation}
where s.l denotes terms which depend on $\Tremb$ for words with a strictly {\it short length} than $l(\bbm )$.\\

The mould $\Tremb$ takes non-trivial values only on resonant words, {\it i.e.} $\bbm \in {\cal A} (F)^*$ such that $\parallel \bbm \parallel.\lambda =0$.
However, the mould $\Semb$ is equal to $1$ on resonant words of length $1$. As a consequence, equation (\ref{trem1}) can not be used to compute the mould $\Tremb$ by
induction on the length of words.\\

Using equation (\ref{etrem2}) we obtain
\begin{equation}
\label{trem2}
\Trem^{\bbm} =\Trem^{\bbm } \Sem^{m_1} \dots \Sem^{m_r}+\mbox{\rm s.l} .
\end{equation}

\subsection{About Ecalle-Vallet results}

All our computations have been done in the alphabet ${\cal D}_{{\cal A} (F)}$. However, J. Ecalle and B. Vallet \cite{ev1} use the initial alphabet
${\cal B}_{A(F)}$ to formulate their results. In order to compare our approach, we first give a simple formula connecting the two alphabets. We then
discuss the essential differences between the moulds $\demb$, $\semb$ and $\tremb$ with the moulds $\Demb$, $\Semb$, and $\Tremb$. The main point is that
contrary to our moulds, Ecalle-Vallet moulds can not be expressed via closed formulae, except for $\demb$.

\subsubsection{Relation between the alphabets ${\cal B}_{A(F)}$ and ${\cal D}_{{\cal A}(F)}$}

By definition, we have the identity
\begin{equation}
1+\di\sum_{n\in A(F)} B_n =\exp \left ( \di\sum_{m\in {\cal A}(F)} D_m \right ) .
\end{equation}
Using the logarithm, we obtain
\begin{equation}
\log \left ( 1+\di\sum_{n\in A(F)} B_n \right ) = \di\sum_{m\in {\cal A}(F)} D_m .
\end{equation}
As $\sum_{n\in A(F)} B_n =\di\sum_{\bn \in A^* (F)} \I^{\bn} B_{\bn}$, we have
\begin{equation}
\di\sum_{\bn \in A^* (F)} (\Log \Ib )^{\bn} B_{\bn} =\di\sum_{m\in {\cal A}(F)} D_m .
\end{equation}
We finally deduce the following relation between ${\cal D}_{{\cal A} (F)}$ and ${\cal B}_{A(F)}$:

\begin{lem}
\label{relationalphabet}
For all $D_m\in {\cal D}_{{\cal A} (F)}$, we have
\begin{equation}
D_m =\di\sum_{\bn \in A(F)^* ,\ \parallel \bn \parallel =m} (\Log \Ib )^{\bn} B_{\bn} .
\end{equation}
\end{lem}

The proof is based on the fact that a differential operator $B_{\bn}$ is of order $\parallel \bn \parallel$.

\subsubsection{The mould $\demb$}

By definition, we have the identity
\begin{equation}
\label{base1}
\di\sum_{\bn \in A(F)^*} \dem^{\bn} B_{\bn} =\di \sum_{m\in {\cal A} (F)\setminus {\cal R}_{{\cal A}(F)}} \di {D_m \over 1-\di e^{m.\lambda}} .
\end{equation}
Using lemma \ref{relationalphabet}, we deduce:

\begin{lem}
The mould $\demb$ of ${\cal M}_{\C} (A(F))$ is defined for all $\bn \in A(F)^*$ by
\begin{equation}
\dem^{\bn} = \di {(-1)^{l(\bn )+1} \over l(\bn )!} \di {1\over 1-\di e^{\parallel \bn\parallel .\lambda }}
\left [ \Ib \right ]^{\bn}_{(\times l(\bn ))} \indic_{N(F)} (\bn ) ,
\end{equation}
where $N(F)=\{ \bn \in A(F)^* ,\parallel \bn \parallel .\lambda \not =0\}$ is the set of non-resonant words of $A(F)^*$ and $\indic_J$ is the
{\it indicatrice} of the set $J$, {\it i.e.} $\indic_J (x)$ is equal to $1$ if $x\in J$, $0$ otherwise.
\end{lem}

This mould is defined directly by Ecalle-Vallet without any details (see \cite{ev1}, p.30).

\begin{proof}
Equation (\ref{base1}) can be rewritten as
\begin{equation}
\label{base2}
\di\sum_{\bn \in A(F)^*} \dem^{\bn} B_{\bn} =\di \sum_{m\in {\cal A} (F)} \di {D_m \over 1-\di e^{m.\lambda}} \indic_{\{m.\lambda \not= 0\}} (m).
\end{equation}
Using lemma \ref{relationalphabet}, we have
\begin{equation}
\label{base3}
\left .
\begin{array}{lll}
\di \sum_{m\in {\cal A} (F)} \di {D_m \over 1-\di e^{m.\lambda}} \indic_{\{m.\lambda \not= 0\}} (m) & = & \di\sum_{m\in {\cal A} (F)}\,
\di\sum_{\bn \in A(F)^* ,\ \parallel \bn \parallel =m} \di { (\Log \Ib )^{\bn} \over 1-\di e^{m.\lambda}} B_{\bn} ,\\
 & = & \di\sum_{\bn \in A(F)^*} \di { (\Log \Ib )^{\bn} \over 1-\di e^{\parallel \bn \parallel .\lambda}}
\indic_{N(F)} B_{\bn} ,
\end{array}
\right .
\end{equation}
using the fact that
\begin{equation}
\bigcup_{m\in {\cal A} (F)} \left \{ \bn \in A(F)^* ,\ \parallel \bn \parallel =m \right \} =A(F)^* ,
\end{equation}
by assumption.

Using lemma \ref{simple} for the mould $\Ib$, we obtain for all $\bn \in A(F)^*$
\begin{equation}
\Log \I^{\bn} =\di {(-1)^{l(\bn )+1} \over l(\bn )!} \left [ \Ib \right ]^{\bn}_{(\times l(\bn ))} .
\end{equation}
Replacing $\Log \Ib$ by its expression in equation (\ref{base3}) we conclude the proof.
\end{proof}

\section{The Poincar\'e-Dulac normal form}

The Trimmed form is constructed using cancellation of non-resonant terms as the classical Poincar\'e-Dulac normal form. However, these two prenormal
forms do not coincide in general. We introduce the universal mould associated to the Poincar\'e-Dulac normal form and the universal mould of the
associated cancellation procedure. The difference between the two procedures lies in the treatment of the homogeneous components of the diffeomorphism.
For a classical approach to the Poincar\'e-Dulac normal form we refer to (\cite{ar} $\S$.B p.178).

\subsection{Homogeneous components and the Trimmed form}

We keep the notations introduced in $\S$.\ref{cancel}. In order to discuss the cancellation of non-resonant terms, we must write our prepared form as
follows:
\begin{equation}
\mbox{\rm Id} +\di\sum_{n\in A(F)} B_n =\exp \mbox{\rm\bf D} =\exp \left ( \di\sum_{m\in {\cal A}(F)} D_m \right ) =\exp \left ( \di\sum_{k\geq 1}
\mbox{\rm\bf D}_k \right ) ,
\end{equation}
where
\begin{equation}
\mbox{\rm\bf D}_k = \di\sum_{n\in {\cal A}(F),\ \mid n\mid =k} D_m ,
\end{equation}
denotes the homogeneous component of degree $k$ of the vector field $\mbox{\rm\bf D}$.\\

For a given vector field $\mbox{\rm\bf D}$ we introduce the following {\it degree of resonance}, denoted by $\mbox{\rm K}$:
\begin{equation}
\mbox{\rm K} =\min_{k\geq 1} \left \{ \mbox{\rm N}_k \not= \emptyset \right \} ,
\end{equation}
where $N_k$ denotes the set of non-resonant letters $m\in {\cal A} (F)$ of degree $k$, {\it i.e.}
\begin{equation}
\mbox{\rm N}_k =\left \{ m\in {\cal A} (F) \, \mid\ \mid m\mid =k ,\ m.\lambda =0 \right \} .
\end{equation}

As a consequence, we have
\begin{equation}
\mbox{\rm\bf D}=\di\sum_{1\leq k <K} \mbox{\rm\bf D}_k +\mbox{\rm\bf D}_K +\di\sum_{k>K} \mbox{\rm\bf D}_k .
\end{equation}
The first sum up to order $K-1$ is made of resonant terms. The first non-resonant terms belong to $\mbox{\rm\bf D}_K$.\\

The field $\mbox{\rm\bf V}$ introduced in $\S$.\ref{sectionsimplified} cancel the non-resonant terms of degree $\mbox{\rm K}$ but introduces several other
terms in the homogeneous components of degree $>\mbox{\rm K}$ which can be non-resonant. As a consequence, even if the field $\mbox{\rm\bf V}$ is constructed in
order to cancel {\it all} the non-resonant terms of the vector field $\mbox{\rm\bf D}$ we have an effective cancellation only for the components of
degree $\mbox{\rm K}$.\\

As a consequence, the vector field $\mbox{\rm\bf V}$ must be modified in order to cancel {\it only} non-resonant terms of degree $\mbox{\rm K}$.

\begin{thm}[Poincar\'e normalization procedure]
\label{tpoincare}
Let $\mbox{\rm\bf S}$ be the vector field defined by the mould
\begin{equation}
\Denb =\left \{
\begin{array}{l}
\di {1 \over 1-\di e^{m .\lambda} } \ \ \mbox{\rm for}\ \ m \in \mbox{\rm N}_{\rm K} (F) ,\\
0\ \ \mbox{\rm otherwise},
\end{array}
\right .
\end{equation}
We denote by $\denb$ the associated mould on ${\cal M}_{\C} (A(F))$, {\it i.e.}
\begin{equation}
\mbox{\rm\bf S} =\di\sumb \Denb \Db =\di\sumb \denb \Bb .
\end{equation}
We call simplified form of $F$ and we denote by $\FPoin$ the automorphism obtained from $F$ under the action of $\exp \mbox{\rm\bf S}$. We have
\begin{equation}
\left .
\begin{array}{lll}
\FPoin & = & \Flin \left ( \di\sum_{\bbm \in {\cal A} (F)^*} \Poin^{\bbm} D_{\bbm} \right ) ,\\
      & = & \Flin \left ( \di\sum_{\bn \in A (F)^*} \poin^{\bn} B_{\bn} \right )
\end{array}
\right .
\end{equation}
with the mould $\Poinb$ given by
\begin{equation}
\Poinb =\edelta \left ( \Exp (\Denb )\right ) \cdot \Exp (\Ib )\cdot \Exp (-\Denb ) ,
\end{equation}
and the mould $\poinb$ given by
\begin{equation}
\poinb =\edelta \left ( \Exp (\denb )\right ) \cdot (\1b +\Ib)\cdot \Exp (-\denb ) .
\end{equation}
\end{thm}

The proof is exactly the same as those of theorem \ref{tsimplifiedform}.

\subsection{The Poincar\'e normal form of order r}

We apply the Poincar\'e normalization procedure inductively in order to cancel non-resonant terms in homogeneous components of higher and higher
degree.

\begin{defi}[Poincar\'e normal form up to order $r$]
Let $r\in \N$, the Poincar\'e normal form up to order $r$ is defined as $\FPoin^r$ obtained from $F$ after $r$ successive simplifications, {\it i.e.}
\begin{equation}
F=\FPoin^0 \stackrel{\AutSimp^1}{\rightarrow} \FPoin^1 \stackrel{\AutSimp^2}{\rightarrow} \dots \stackrel{\AutSimp^r}{\rightarrow} \FPoin^r ,
\end{equation}
where $\AutSimp^i$ is the automorphism of simplification defined by
\begin{equation}
\AutSimp^i =\exp ( \mbox{\rm\bf S}_i ) ,
\end{equation}
with $\mbox{\rm\bf S}_i$ the vector fields associated to the mould $\Denb$ on the alphabet ${\cal A} (\FPoin^{i-1} )$ associated to $\FPoin^{i-1}$.
\end{defi}

Using theorem \ref{tpoincare}, we obtain:

\begin{thm}
For all $r\in \N$, the Poincar\'e normal form up to order $r$ denoted $\FPoin^r$ possesses a mould expansion, {\it i.e.} there exist moulds denoted by
${}_r \Poinb \in {\cal M}_{\C} ({\cal A}(F))$ and ${}_r \poinb \in {\cal M}_{\C} (A(F))$ such that
\begin{equation}
\FPoin^r =\Flin \left ( \di\sumb {}_r \Poinb \Db \right ) =\Flin \left ( \di \sumb {}_r \poinb \Bb \right ) .
\end{equation}
\end{thm}

As for the moulds ${}_r \semb$ and ${}_r \Semb$, we have explicit inductive formulae to compute the moulds ${}_r \poinb$ and ${}_r \Poinb$ using
only $\poinb$ and $\Poinb$.

\subsection{The Poincar\'e-Dulac normal form}

The mould formulation of the Poincar\'e-Dulac normal form is:

\begin{defi}
The Poincar\'e-Dulac normal form of $F$ is the limit of the Poincar\'e normalization procedure.
\end{defi}

\begin{thm}
The Poincar\'e-Dulac normal form is a continuous prenormal form given by
\begin{equation}
\left .
\begin{array}{lll}
\FDulac & = & \Flin \left ( \di\sum_{\bbm \in {\cal A} (F)^*} \Dulac^{\bbm} D_{\bbm} \right ),\\
 & = & \Flin \left ( \di\sum_{\bn \in A (F)^*} \dulac^{\bn} B_{\bn} \right )
\end{array}
\right .
\end{equation}
with the moulds $\Dulacb$ and $\dulacb$ defined by
\begin{equation}
\begin{array}{c}
\Dulacb -\1b =\limstat_{r\rightarrow \infty} \left  [ \Poinb -\1b \right ] ^{(\circ r)} ,\\
\dulacb -\1b =\limstat_{r\rightarrow \infty} \left  [ \poinb -\1b \right ] ^{(\circ r)} ,
\end{array}
\end{equation}
where $\limstat$ is the stationary limit.
\end{thm}

The mould $\Dulacb$ (or $\dulacb$) is the {\it universal} part of the Poincar\'e-Dulac normal form as it does not depends on the exact values of the coefficients
coming in the Taylor expansion of the diffeomorphism. It seems impossible to characterize such kind of object without using moulds.




\begin{appendix}
\section{About the Baker-Campbell-Hausdorff formula}

The Baker-Campbell-Hausdorff formula covers at least two formulae which are of interest for the computation of
continuous prenormal forms for vector fields and diffeomorphisms. \\

Let $A$ and $B$ be two linear operators. We denote by
$$\exp A=\di\sum_{k\geq 0} \di {A^k \over k!} ,$$
where $A^k=A\circ \dots \circ A$, $k$ times.
The Baker-Campbell-Hausdorff formula is given by
$$
(\exp A )\, . B\, . (\exp (-A)) =\di\sum_{m\geq 0} \di {B_m \over m!} ,
\eqno{(BCH_1 )}$$
where
\begin{equation}
B_m =[A,B]_m =[A[A,\dots ,[A,B]]\dots ] ,
\end{equation}
with the convention that $B_0 =B$.

A consequence of this formula is
$$ \exp A \exp B =\exp (A\star B ),\eqno{(BCH_2)}$$
where
\begin{equation}
A\star B =A+B+\di {1\over 2} [A,B]+\di{1\over 12} [A,[A,B]] -{1\over 12}[B,[A,B]] +\dots
\end{equation}

\section{Proof of theorem \ref{explicitSem}}
\label{proofSem}

In order to compute the mould $\Semb$, we first compute $\Exp \Ib \cdot \Exp (-\Demb )$. We have

\begin{equation}
\left .
\begin{array}{lll}
\left ( \Exp \Ib \cdot \Exp (-\Demb ) \right )^{\bn} & = & \di\sum_{\bn^1 \bn^2 =\bn} (\Exp \Ib)^{\bn^1} \Exp (-\Demb )^{\bn^2} ,\\
   & = & \di\sum_{\bn^1 \bn^2 =\bn} \left ( \un^{\bn^1} +\di {1\over l(\bn^1 )!} \left [\Ib \right ]_{(\times l(\bn^1 ))}^{\bn^1} \right ) \left ( \un^{\bn^2} +\di {(-1)^{l(\bn^2)}\over l(\bn^2 )!} \left [\Demb \right ]_{(\times l(\bn^2 ))}^{\bn^2} \right ) ,\\
   & = & \di\sum_{\bn^1 \bn^2 =\bn} \bigg ( \un^{\bn^1} \un^{\bn^2} + \un^{\bn^1} \di {(-1)^{l(\bn^2)}\over l(\bn^2 )!}
   \left [\Demb \right ]_{(\times l(\bn^2 ))}^{\bn^2}  \\
    & & + \un^{\bn^2} \di {1\over l(\bn^1 )!} \left [\Ib \right ]_{(\times l(\bn^1 ))}^{\bn^1} + \di {(-1)^{l(\bn^2)}\over l(\bn^1)! l(\bn^2 )!} \left [ \Ib \right]_{(\times l(\bn^1))}^{\bn^1} \left [\Demb \right ]_{(\times l(\bn^2 ))}^{\bn^2} \bigg ) .
\end{array}
\right .
\end{equation}

It is clear that~$\left ( \Exp \Ib \cdot \Exp (-\Demb ) \right )^\emptyset = 1$. If~$l(\bn) \ge 1$ we have

\begin{equation}
\left .
\begin{array}{lll}
\left ( \Exp \Ib \cdot \Exp (-\Demb ) \right )^{\bn} & = & \di {(-1)^{l(\bn)}\over l(\bn )!} \left [\Demb \right ]_{(\times l(\bn ))}^{\bn} + \di {1\over l(\bn )!} \left [\Ib \right ]_{(\times l(\bn ))}^{\bn} \\
    & &  \quad +\di\sum_{\bn^1 \bn^2 =\bn \atop \bn^1 \ne \emptyset} \left ( {(-1)^{l(\bn^2)} \over l(\bn^1)! l(\bn^2 )!} \left [\Demb \right ]_{(\times l(\bn^2 ))}^{\bn^2} \right) ,\\
    & = & \di {(-1)^{l(\bn)}\over l(\bn )!} \left [\Demb \right ]_{(\times l(\bn ))}^{\bn} + \di {1\over l(\bn )!} \\
    & & + \di\sum_{j=d(\bn )+1}^{l(\bn )+1}
\di { (-1)^{l(\bn^{\ge j})} \left [\Demb \right ]^{\bn^{\geq j}}_{(\times l(\bn^{\geq j} ))} \over l(\bn^{<j} )! l(\bn^{\geq j} )! } .
\end{array}
\right .
\end{equation}

Now we can compute~$\Semb$.

\begin{equation}
\left .
\begin{array}{lll}
\Sem^{\bn} & = & \left ( \edelta \left ( \Exp (\Demb )\right ) \cdot \Exp (\Ib )\cdot \Exp (-\Demb ) \right )^{\bn} ,\\
    & = & \di \sum_{\bn^1 \bn^2 =\bn} \left ( \edelta  \Exp (\Demb ) \right )^{\bn^1} \bigg ( \di {(-1)^{l(\bn^2)}\over l(\bn^2 )!} \left [\Demb \right ]_{(\times l(\bn^2 ))}^{\bn^2} \\
    &  & \quad \quad \quad \quad \quad + \di {1\over l(\bn^2 )!} + \di\sum_{j=d(\bn^2 )+1}^{l(\bn^2 )+1}
\di { (-1)^{l((\bn^2)^{\ge j})}\left [\Demb \right ]^{(\bn^2)^{\geq j}}_{(\times l((\bn^2)^{\geq j} ))} \over l((\bn^2)^{<j} )! l((\bn^2)^{\geq j} )! } \bigg ) ,\\
    & = &  \di \sum_{\bn^1 \bn^2 =\bn} e^{-\lambda. \parallel \bn^1 \parallel}\left ( 1^{\bn^1} + {1 \over l(\bn^1)!} \left [\Demb \right ]_{(\times l(\bn^1 ))}^{\bn^1} \right ) \times \\
    & & \quad \bigg ( \di {(-1)^{l(\bn^2)}\over l(\bn^2 )!} \left [\Demb \right ]_{(\times l(\bn^2 ))}^{\bn^2} + \di {1\over l(\bn^2 )!} + \di\sum_{j=d(\bn^2 )+1}^{l(\bn^2 )+1}
\di { (-1)^{l((\bn^2)^{\ge j})}\left [\Demb \right ]^{(\bn^2)^{\geq j}}_{(\times l((\bn^2)^{\geq j} ))} \over l((\bn^2)^{<j} )! l((\bn^2)^{\geq j} )! }  \bigg ) ,\\
    & = &  \di {(-1)^{l(\bn)}\over l(\bn )!} \left [\Demb \right ]_{(\times l(\bn ))}^{\bn} + \di {1\over l(\bn )!} + \di\sum_{j=d(\bn )+1}^{l(\bn )+1}
\di { (-1)^{l(\bn^{\ge j})} \left [\Demb \right ]^{\bn^{\geq j}}_{(\times l(\bn^{\geq j} ))} \over l(\bn^{<j} )! l(\bn^{\geq j} )! }  \\
    & & + \di \, e^{-\lambda. \parallel \bn \parallel} \, \mbox{\rm 1}^{\bn}
  +\di\sum_{i=1}^{q(\bn ) \wedge (l(\bn)-1)} \di {e^{-\lambda. \parallel \bn^{\leq i} \parallel } \over l(\bn^{\leq i} )!} \left [\Demb \right ]
^{\bn^{\leq i}}_{(\times l(\bn^{\leq i} ))} \times \\
    & & \left (
 \di  {(-1)^{l(\bn^{>i})} \over l(\bn^{>i} )!} \left [ \Demb \right ]^{\bn^{>i}}_{(\times l(\bn^{>i} ))} +\di {1\over l(\bn^{>i} )!} +
 \di\sum_{j=d(\bn^{>i} )+1}^{l(\bn^{>i} )+1}
 \di {(-1)^{l(\bn^{>i})^{\geq j}} \left [\Demb \right ] ^{(\bn^{>i} )^{\geq j}}_{(\times l((\bn^{>i} )^{\geq j} ))} \over l((\bn^{>i} )^{<j} )! l((\bn^{>i})^{\geq j} )!} \right ) .
\end{array}
\right .
\end{equation}
This concludes the proof.

\end{appendix}

\vskip 1cm
\parindent 0pt
\begin{tiny}
Jacky CRESSON
\vskip 2mm
Universit\'e de Pau et des Pays de l'Adour,

Laboratoire de Math\'ematiques appliqu\'ees de Pau, CNRS UMR 5142

jacky.cresson@univ-pau.fr
\vskip 2mm
{\it and}
\vskip 2mm
Institut des Hautes \'Etudes Scientifiques (I.H.\'E.S)

Le Bois-Marie, 35 Route de Chartres, F-91440 Bures sur Yvette, France.

cresson@ihes.fr
\vskip 7mm

Jasmin RAISSY

University of Pisa

raissy@mail.dm.unipi.it
\end{tiny}
\end{document}